\documentclass[11pt]{article}
\newtheorem{theorem}{Theorem}

\newtheorem{claim}{Claim}

\newtheorem{lemma}{Lemma}

\usepackage{amsmath,amssymb,color,graphicx,epsf,verbatim}

\usepackage[left=0.99in, right=0.99in, top=1in, bottom=0.99in]{geometry}

\def\reals{{\mathbb R}}

\def\proof{\noindent{\em Proof: }}
\def\qed{$\spadesuit$}

\begin{document}
\begin{titlepage}
\title{Beyond the Richter-Thomassen Conjecture}

\author{J\'anos Pach\thanks{EPFL, Lausanne and R\'enyi Institute,
    Budapest. Supported by Swiss National
    Science Foundation Grants 200020-144531 and 20021-137574. Email: {\tt pach@cims.nyu.edu}}
  \and Natan Rubin\thanks{Ben Gurion University of the Negev,  Beer-Sheba, Israel.
Email: {\tt rubinnat.ac@gmail.com}. Supported by Minerva Fellowship Program of the Max Planck Society, by the Fondation Sciences Math\'{e}matiques de Paris (FSMP), and by a public grant overseen by the French National Research Agency (ANR) as part of the "Investissements d'Avenir" program (reference: ANR-10-LABX-0098).}
  \and G\'abor Tardos\thanks{R\'enyi Institute, Budapest. Supported by the ``Lend\"ulet'' Project of the Hungarian Academy of Sciences. Email: {\tt
      tardos@renyi.hu}}}

\maketitle
\begin{abstract}
If two closed Jordan curves in the plane have precisely one point in common, then it is called a {\em touching point}. All other intersection points are called {\em crossing points}. The main result of this paper is a Crossing Lemma for closed curves: In any family of $n$ pairwise intersecting simple closed curves in the plane, no three of which pass through the same point, the number of crossing points exceeds the number of touching points by a factor of at least $\Omega((\log\log n)^{1/8})$.

\smallskip

As a corollary, we prove the following long-standing conjecture of Richter and Thomassen: The total number of intersection points between any $n$ pairwise intersecting simple closed curves in the plane, no three of which pass through the same point, is at least $(1-o(1))n^2$.
\end{abstract}

\maketitle
\end{titlepage}
\section{Introduction}
\noindent{\bf Arrangements of curves and surfaces.} It was a fruitful and surprising discovery made in the 1980s that the Piano Mover's Problem and many other algorithmic and optimization questions in motion planning, ray shooting, computer graphics etc., boil down to computing certain elementary substructures (e.g., cells, envelopes, $k$-levels, or zones) in arrangements of curves in the plane and surfaces in higher dimensions~\cite{Ed87,KLPS86,PaS09,ShA95}. Hence, the performance of the most efficient algorithms for the solution of such problems is typically determined by the {\em combinatorial complexity} of a single cell or a collection of several cells in the underlying arrangement, that is, the total number of their vertices, edges, and faces of all dimensions. 

The study of arrangements has brought about a renaissance of Erd\H os-type combinatorial geometry.
For instance, in the plane, Erd\H os's famous question~\cite{Er46} on the maximum number of times the unit distance can occur among $n$ points in the plane can be generalized as follows~\cite{CEGSW90}: {\it What is the maximum total number of sides of $n$ cells in an arrangement of $n$ unit circles in the plane?} In the limiting case, when $k$ circles pass through the same point $p$ (which is, therefore, at unit distance from $k$ circle centers), $p$ can be regarded as a degenerate cell with $k$ sides.


Several beautiful paradigms have emerged as a result of this interplay between combinatorial and computational geometry, from the random sampling argument of Clarkson and Shor~\cite{CS89} through epsilon-nets (Haussler-Welzl~\cite{HW87}) to the discrepancy method (Chazelle~\cite{Cha00}).
It is worth noting that most of these tools are restricted to families of curves and surfaces of {\it bounded description complexity}, in the sense that the arrangement of any constantly many objects is composed of constantly many cells.

\smallskip

Another tool that proved to be applicable to Erd\H os's questions on repeated distances is Sz\'ekely's {\em Crossing Lemma}~\cite{Sz97}. It states that no matter how we a draw a sufficiently large graph $G=(V,E)$ in the plane or on a fixed surface, the number of crossings between its edges is at least $\Omega(|E|^3/|V|^2)$. In particular, this implies that if $G$ has a lot more edges than vertices, then its number of crossings is much larger than its number of edges. The best known results on the $k$-set problem~\cite{De98}, needed for the analysis of many important geometric algorithms, and the most elegant proofs of the Szemer\'edi-Trotter theorem~\cite{SzT83a}, \cite{SzT83b} on the maximum number of incidences between a set of points and a set of lines (or other, more complicated, curves) were also established using the Crossing Lemma~\cite{PaS98}. These proofs easily generalize from lines to pseudo-segments (i.e., curves with at most one intersection per pair).

\medskip

\noindent{\bf Tangencies and lenses.} Motivated by potential applications in motion planning, Tamaki and Tokuyama~\cite{TT98} extended the $k$-set bounds and incidence bounds from lines to more general curves, by trying to cut the curves into as few pseudo-segments as possible, and then applying the known bounds to them. In this context, the number of {\em tangencies (touchings)} between the original curves plays a special role. By locally perturbing two curves in a small neighborhood of their touching point, one can create two nearby crossings and a small ``lens'' between them. In order to decompose the curves into pseudo-segments, we have to make at least one cut on the boundary of each lens. In many scenarios, the number of cuts needed is roughly proportional to the number of touching points, more precisely, to the maximum number of non-overlapping lenses. This approach was later refined and extended in a series of papers~\cite{Chan1}, \cite{Chan2}, \cite{Chan3}, \cite{ArS02}, \cite{AgS05}, \cite{MaT06} and \cite{PseudoCircles}.

In particular,  Agarwal {\it et al.} \cite{PseudoCircles} studied arrangements of pseudo-discs (that is, closed Jordan curves with at most two intersections per pair) and used lenses to establish several fundamental results on geometric incidences and 
cell complexity.
Their analysis crucially relied on the following claim: Any family of $n$ pairwise intersecting pseudo-circles admits at most $O(n)$ tangencies. In the special case where the curves are algebraic, any incidence or tangency can be described by a polynomial equation. 
Following the pioneering work of Dvir~\cite{Dv10}, Guth and Katz~\cite{GK10}, \cite{GK15}, many of these problems have been revisited from an algebraic perspective.

The structure of tangencies between {\em convex sets} was addressed in \cite{PST12}. It was shown that the number of tangencies between $n$ members of any family of plane convex sets that can be obtained as the union of $k$ packings (systems of disjoints sets) is at most $O(kn)$. The proof of this fact is somewhat delicate, because the boundaries of two convex sets can cross any number of times.

\medskip
\noindent{\bf Our result.} The main result of this paper is a Crossing Lemma for a family of $n$ pairwise intersecting closed curves. We are going to show, roughly speaking, that the number of {\it proper} (i.e., transversal) crossings between the curves is much larger than the number of touching pairs of curves, provided that $n$ is sufficiently large.

To formulate this result more conveniently, we need to agree on the terminology. We say that two (open or closed) curves {\em intersect} if they have at least one point in common. An intersection point $p$ is called a {\em touching point} (in short, a {\em touching}) if $p$ is the {\em only} intersection point of the two curves, and they do not properly cross at $p$. All other intersection points will be referred to as {\em crossing points} (in short, {\em crossings}). Note that this definition is somewhat counterintuitive: if two curves have two points of tangencies, we call both of them a crossing.\footnote{Alternatively, we can make the pair touch at a single point of tangency by slightly perturbing them in the neighborhood of the other touching.} It is assumed throughout that all curves are in {\em general position}, that is, no three of them pass through the same point and no two share infinitely many points.

Now we can state our Crossing Lemma for closed Jordan curves:

\begin{theorem}\label{main} 
Let $A$ be a collection of $n$ pairwise intersecting closed Jordan curves in the plane, no three of which pass through the same point. Let $T$ denote the set of touching points and let $X$ denote the set of crossing points between the elements of $A$. We have
\begin{equation}\label{Eq:Gap}
|X|=\Omega\left(|T|(\log\log n)^{1/8}\right).
\end{equation}
\end{theorem}

At first glance, one might believe that the statement remains true even if we drop the assumption that the curves are pairwise intersecting. The following example will convince us that this is not the case. Take $n-k$ pairwise disjoint circles in the plane. It is easy to select $k$ other closed curves in general position such that each of them touches every circle and any pair of them cross at most $n-k$ times. In this arrangement, $|T|=k(n-k)$ and $|X|\le{k\choose 2}(n-k)$, so that we have $\frac{|X|}{|T|}\le\frac{k-1}{2},$ which does not necessarily tend to infinity as $n$ increases.

\smallskip
We use Theorem~\ref{main} to prove the following long-standing conjecture of Richter and Thomassen~\cite{RiT95}: \\
 {\em Given a collection of $n$ pairwise intersecting closed curves in general position in the plane, $|X|+|T|$, the total number of intersection points between the curves, is at least $(1-o(1))n^2$}. Note that if there are no touchings between the curves, then any two curves cross at least twice, so that the number of intersection points is at least $2{n\choose 2}=(1-o(1))n^2$. However, if touchings are allowed, the situation is more complicated.

The best known general lower bound, $|X|+|T|\ge(4/5-o(1))n^2$ was shown by Mubayi~\cite{Mu02}. The Richter-Thomassen conjecture was confirmed by Salazar \cite{Sa99} in the case when any pair of curves have at most a bounded number of points in common. In an earlier paper~\cite{PRT15}, the authors settled the special case where the curves are {\em convex} or, more generally, if each curve can be cut into a constant number of $x$-monotone arcs. (An arc is called {\em $x$-monotone} if every vertical line intersects it in at most one point.) The problem has remained open for general families of simple closed curves. 

Now we can prove the conjecture in the general case, by a direct application of Theorem 1. 

\begin{theorem}\label{richterthomassen}
The total number of intersection points between $n$ pairwise intersecting closed curves in general position in the plane is at least $(1-o(1))n^2$.
\end{theorem}

To deduce Theorem~\ref{richterthomassen} from our new Crossing Lemma for curves (Theorem~\ref{main}), it is enough to notice that if $|T|=o(n^2)$, then the statement follows from the trivial bound $|X|\ge2({n\choose2}-|T|)$. Otherwise, if $|T|\ge \varepsilon n^2$ for some $\varepsilon>0$, Theorem~\ref{main} immediately implies that $|X|\ge\varepsilon n^2(\log\log n)^{1/8},$ which is much better than required.

\smallskip
\noindent{\bf A remark and open problems.} Instead of {\it closed} Jordan curves, Theorem 1 remains valid for any family of pairwise intersecting (open) {\it Jordan arcs}. To see this, replace each arc by a slightly inflated closed curve so that during the process all tangencies are preserved. The number of crossings between the new curves may increase by a factor of at most $4$. Applying Theorem 1 to this new family, we obtain that the ratio between the number of crossings and the number of touchings determined by the original arcs is at least $\Omega\left((\log\log n)^{1/8}\right)$.  

Fox {\em et al.} \cite{FFPP10} constructed two $n$-size families $A$ and $B$ of pairwise intersecting $x$-monotone arcs in the plane such that every curve in $A$ touches every curve in $B$, and the total number of crossings between the members of $A\cup B$ is $O(n^2\log n)$. They showed that in this setting the bound is optimal. We conjecture that our Theorem 1 can also be improved as follows.

\medskip
\noindent{\bf Conjecture 3.} {\em Let $A$ be a collection of $n$ pairwise intersecting Jordan arcs in the plane, no three of which pass through the same point. Let $T$ denote the set of touching points and $X$ the set of crossing points between the elements of $A$. We have
$|X|=\Omega(|T|\log n)$.}

The proof of Theorem 1 makes substantial use of the property that any two curves intersect. In the case of $x$-monotone arcs, it was shown in \cite{PRT15} that a similar crossing lemma holds under the weaker assumption that there are many intersecting pairs of arcs. We believe that the same is true for closed Jordan curves.     

\medskip
\noindent{\bf Conjecture 4.} {\em Let $A$ be a collection of $n$ closed Jordan curves in the plane, no three of which pass through the same point. Let $T$ denote the set of touching points and $X$ the set of crossing points between the elements of $A$. There is a function $f(x)$ with $\lim_{x\rightarrow\infty}f(x)=\infty$ such that $|X|=\Omega(|T|f(|T|/n))$.}

\smallskip
Perhaps this last conjecture holds even for $f(x)=\log x$.  The first step towards its resolution could be to prove Conjecture 4 for collections of curves with at least $\varepsilon n^2$ intersecting pairs.

\medskip
\noindent{\bf Algebraic techniques.}
As mentioned before, the polynomial technique of Guth and Katz \cite{GK10,G2}, which led to a spectacular breakthrough concerning Erd\H os's problem on distinct distances, has inspired a lot of recent research related to incidences between points, curves, and surfaces \cite{G3,G1,GK15,SSZ15}. Unfortunately, the new techniques only apply in an algebraic framework, where the curves and surfaces in question must be algebraic varieties of bounded degree. Since two algebraic curves of bounded degree that do not share a component have only a bounded number of points in common, restricting the Richter-Thomassen conjecture to such curves, reduces the question to a relatively simple combinatorial problem. This was pointed out by Salazar [Sa99]. For many similar problems related to intersection patterns of curves, including the Erd\H os-Hajnal conjecture [EH] for intersection graphs of curves in the plane, our present techniques are not sufficient to handle the case when two curves may intersect an arbitrary number of times [FPT11, FP08, FP10]. There are only very few exceptional examples, when one is able to drop this assumption [Ma14, FP12]. The main result of this paper represents one of these rare exceptions.

\medskip

\noindent{\bf Overview of the proof.} The proof of our Crossing Lemma (Theorem 1) requires a delicate {\em charging scheme}: a weighted assignment of of the touchings (elements of $T$) to the crossings (elements of $X$) between the curves. In order to prove (1), average amount of charge received by a touching $t\in T$ from its neighboring crossings has to be much larger than the average amount of charge ``paid'' by a crossing. To estimate a charge paid by a crossing $x$ between two curves $a,b\in A,$ we use the local topology of of the underlying arrangement, and we associate $x$ with a ``lens'' formed by $a$ and $b$. Our proof was inspired by charging strategies used for estimating the combinatorial complexity of various special features of bounded ``description complexity'' in various geometric arrangements; see [ShA95] for a comprehensive survey.

\section{Proof of Theorem \ref{main}}

Let $\alpha_1=(\log\log n)^{1/8}/10$, where $\log$ denotes the binary
logarithm. We need to prove $|X|/|T|=\Omega(\alpha_1)$. For this proof we
assume that $|X|\le\alpha_1{n\choose2}$ as otherwise $|X|>\alpha_1|T|$
trivially holds. Note that as the statement we want to prove is asymptotic we
can simplify our calculations by assuming that $n$ is large enough.

\begin{figure}[htbp]
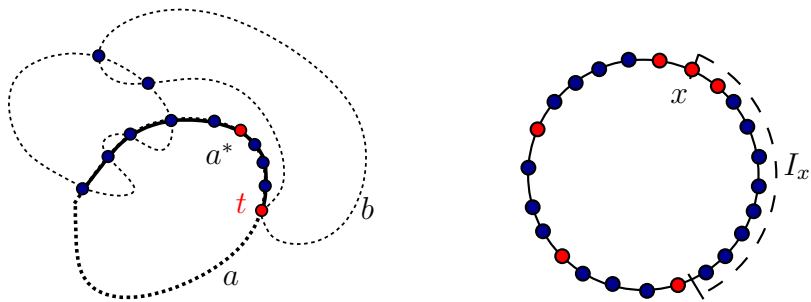

\begin{center}
\input{HappyTouching.pstex_t}\hspace{2cm}\input{Covering.pstex_t}
\caption{\small Left: The touching $t$ is $a$-happy. The points of $T$ (resp.,
  $X$) along $a$ are colored red (resp., blue). The arc $a^*\subset a$
  satisfies $|X\cap a^*|\geq \alpha_1 |T\cap a^*|$. Right: Lemma
  \ref{Lemma:RedBlue}. Each red point $x\in R$ is an endpoint of an arc $I_x$
  that satisfies $w(B\cap I_x)\geq \lambda |R\cap I_x|$. (In the depicted
  scenario $w(x)=1$ for all $x$ and we have $\lambda=4$.)}
\label{Fig:Covering}
\end{center}
\end{figure}

We say that a touching point $t\in T$  is {\em $a$-happy} if  $a$ is one the
two curves in $A$ containing $t$ and $t$ is an end point of a closed arc $a^*\subset a$ with $|X\cap a^*|\geq
\alpha_1|T\cap a^*|$. See Figure \ref{Fig:Covering} (left). We say that $t$ is {\em happy} if it is $a$-happy for some
$a\in A$ with $t\in a$. We call a touching point {\em sad} if it
is not happy and denote the set of sad touching points by $T'$.
We bound the
number of happy and sad touching points (relative to the number of crossing
points) separately.

In an appropriate size ``neighborhood'' of a happy touching point the
crossing-to-touching ratio is high, so bounding the number of happy touching
points is relatively simple. We will use the following lemma. (See Figure \ref{Fig:Covering} (right). Notice that the sets $R$ and $B$ may
overlap.)

\begin{lemma}\label{Lemma:RedBlue}
Let $a$ be a simple (open or closed) Jordan curve and let $R,B\subset
a$ be two finite subsets of points on $a$. Let $\lambda$ be a positive
constant and let $w: B\rightarrow \reals$ be a
positive weight function. For $S\subseteq B$ the weight of $S$ is $w(S)=\sum_{x\in S}w(x)$. If every point $x\in R$ is the
endpoint of a closed 
arc $I_x\subseteq a$ that satisfies $w(B\cap I_x)\geq \lambda |R\cap I_x|$,
then we have $w(B)\geq \lambda|R|/3$.
\end{lemma}

\medskip
\proof We prove by induction on $|R|$. The claim trivially holds if $R$
is empty, so we assume $R$ is not empty and the statement of the lemma holds
for $R'$ and $B'$ as long as $|R'|<|R|$.

Let us choose $x\in R$ to maximize $|R\cap I_x|$ breaking ties arbitrary. Let
$B'=B\setminus I_x$ and $R'=\{y\in R\mid I_y\cap I_x=\emptyset\}$. For $y\in
R'$ we have $B'\cap I_y=B\cap I_y$ and $R'\cap I_y\subseteq R\cap I_y$, so the
assumption of the lemma is satisfied for $R'$ and $B'$. As $x\notin R'$ we have
$|R'|<|R|$ and thus, by the inductive hypothesis, we have
$w(B')\ge\lambda|R'|/3$. By the choice of $x$ every  $y\in R\setminus R'$
must either be in $I_x$ or it is one of the $|I_x\cap R|$ next points in $R$
in either side of the arc $I_x$. So we have $|R|-|R'|\le3|I_x\cap R|$. We
further have $w(B)-w(B')=w(I_x\cap B)\ge\lambda|I_x\cap
R|\ge\lambda(|R|-|R'|)/3$. Adding this inequality to the one obtained from the
inductive hypothesis finishes the proof. \qed
\medskip

\begin{lemma}\label{happy}
$|T|-|T'|\le 6|X|/\alpha_1$.
\end{lemma}
\medskip

\proof We apply Lemma \ref{Lemma:RedBlue} for each curve $a\in A$ with
$\lambda=\alpha_1$, $B=B_a=a\cap X$ and $R=R_a$ being the set of $a$-happy
touching points. We use the uniform weight function $w(x)=1$ for each $x\in
B_a$. We obtain $|B_a|\ge\alpha_1|R_a|/3$. Summing this for all
$a\in A$ we get $2|X|$ on the left hand side and at least $\alpha_1(|T|-|T'|)/3$
on the right hand side. \qed

\medskip
The overall setup for the bound on the sad touching points is the so called charging method: we send
certain amounts of ``charge'' from points in $X$ to points in $T'$. If we
manage to make sure that the total charge sent by any point in $X$ is at most
$c_{\mbox{\scriptsize out}}$ and the total charge received by any point of $T'$ is
at least $c_{\mbox{\scriptsize in}}$, then we have established that $|X|/|T'|\ge
c_{\,mbos{\scriptsize in}}/c_{\mbox{\scriptsize out}}$. Note that we used the same method in our paper \cite{PRT15}
to prove certain special cases of the Richter--Thomassen conjecture.

Our charging is done in phases: in each phase we fix the value of the parameter $k$ (``the
scale'') and perform certain chargings with that scale. Our goal is to make sure that each point
in $X$ sends out a constant charge in each phase, while each touching in $T'$
receives a charge of $\Omega(\alpha_1)$ in each phase. If we could do this, then a single phase would be enough
to prove Theorem~\ref{main}. But we will not quite achieve this goal. The first minor technical difficulty is
that some of our charging rules sends charge not only from points in $X$
but also from points in $T$. A more important problem is  that some points in $X$ are overcharged in certain
rounds: they send out more than a constant amount of charge. This problem is solved by considering
several rounds at once and showing that {\em on average} no crossing is overcharged.

We introduce some notation. We orient each curve $a$ in $A$ so
that all other curves from $A$ touching $a$ touches it on its right side. This
is possible as if $a$ has a touching curve on either side, then these curves
are not intersecting counter to our assumption that each pair of curves in
$A$ intersect. We use the word {\em arc} for closed segments of the curves in $A$. We will use lowercase
letters with an asterisk to denote arcs. The arcs inherit
their orientation from the curve of $A$ containing it and this orientation distinguishes the
{\em starting point} and the {\em end point} of an arc. For distinct points $p$ and $q$ in a curve
$a\in A$ we write {\em the arc of $a$ from $p$ to $q$} to refer to the single arc on $a$ with $p$
as its starting and $q$ as its end point. We can simply refer to an arc as ``the arc from $p$ to $q$''
unless $p,q\in X$ represent two intersections of the same two curves from $A$.
The orientation also makes references like ``the next $k$
points of $T$ along $a$ after $p$'', or ``the last $k$ points of $T$ along $a$ before $p$'' 
unambiguous. By the {\em length} of an arc
we mean the number of sad touching points it contains. Let $x\in X$ be an intersection between the curves
$a,b\in A$ and let $y$ be the another intersection point of the same two curves, the next such point along
$a$. We call the arc of $a$ from $x$ to $y$ a {\em lens}. (Note that most texts include
both arcs from $x$ to $y$ in their definition of a lens, but for us it is
simpler to focus on a single arc.) We write $X'$ for all the intersection
points of the curves in $A$: $X'=X\cup T$.

We set the following parameters: $\alpha=\alpha_1+2$, $v=21000\alpha^8$
depending only on $n$ and the parameter $w=w(k)=k^3/(2000\alpha^5n^2)$ that
also depends on the scale.

\begin{figure}[htbp]
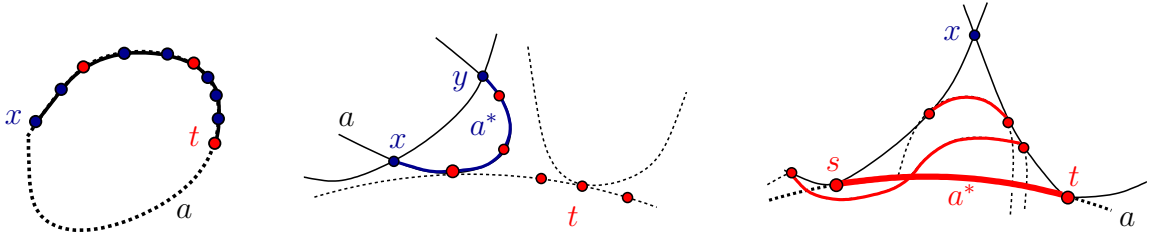

\begin{center}
\input{FirstRule.pstex_t}\hspace{1cm}\input{ChargingRuleLens.pstex_t}\hspace{1cm}\input{LastRuleApex.pstex_t}
\caption{\small Left: First charging rule. A point $x\in X'$ pays $1/k$ units
  to a sad touching $t\in T'$ if the interval from $t$ to $x$ (or vice versa)
  has length at most $k$. Center: Second charging rule. The touching $t$
  receives $v/(k(l+w))$ units from the lens $a^*$ with endpoints $x$ and
  $y$. Right: Third charging rule. The endpoint $t$ of the poor arc $a^*$
  receives $2a/k$ units from the apex $x\in X'$ because there exist at most $k/a$ other poor arcs with apex $x$ that start between $x$ and $t$.}
\label{Fig:Rules}
\end{center}
\vspace{-0.5cm}
\end{figure}

Let us consider the phase with scale $k$. We start with describing our three
charging rules sending charges from intersection points in $X'$ to sad
touchings in $T'$. See Figure \ref{Fig:Rules}.

\smallskip
{\bf First charging rule} A point $x\in  X'$ sends a charge of $1/k$ to a
point $t\in T'$ if the two points are on a common curve of $A$ and
either the arc from $x$ to $p$ or the arc from $p$ to $x$ has length at most $k$.

{\bf Second charging rule} If the length $l$ of a lens $a^*$ satisfies
$l\le3\alpha^3k$, then $a^*$ sends a charge of $v/(k(l+w))$ to all
points
$t\in T'$ that have an arc of length at most $k+1$ from $t$ to a point in
$a^*\cap T'$ and this arc is not along the same curve as $a^*$.

\smallskip
For accounting purposes we consider a charge sent by a lens $a^*$ to be sent
by the starting point of $a^*$.

We call a point of $T'$ {\em poor} in this phase if it receives less than a
total charge of $\alpha$ from the first two charging rules. We call an arc
{\em poor} if it starts at a poor touching point,
ends at a sad touching point and has length at most $k+1$.

Let $a^*$ be an arc of a curve $a\in
A$ starting at $s\in T'$ and ending at $t\in T'$. Let $b$ and $b'$ be the curves in $A$ touching
$a$ in the points $s$ and $t$, respectively. We define the {\em apex} of the arc
$a^*$ as the first point on $b'$ after $t$ that also belongs to $b$. This is a
well defined point in $X'$ as $b$ and $b'$ (as any pair of curves in $A$) must
intersect.

\smallskip
{\bf Third charging rule} Let $a^*$ be a poor arc starting at $t\in T'$ and
having $x\in X'$ as its apex. The intersection point $x$ sends a charge of $2\alpha/k$ to $t$ in this phase unless there are
more than $k/\alpha$ poor arcs, each starting at a point in the arc from $x$
to $t$ and having $x$ as its apex.
\bigskip

\vspace{-0.3cm}
\subsection{Total charge sent}

\begin{lemma}\label{out13}
The total charge sent from a intersection point $x\in X'$ in a phase according
to the first and third rules is at most $8$.
\end{lemma}

\proof The first rule sends a charge of $1/k$ from $x$ to the first $k$ sad
touching points in each of four ``directions'' (in both directions of both
curves containing $x$). That is at most $4k$ sad touching points
for a total charge of at most $4$.

The third rule sends a charge of $2\alpha/k$ to the first $\lfloor k/\alpha\rfloor$ touching points from $x$
along either curves containing $x$ satisfying a certain condition (namely being the starting point of
a poor arc having $x$ as its apex) for a total of at most $4$.
\qed

A statement similar to Lemma~\ref{out13} is false for the second charging rule because it severely overcharges
the lenses whose length is approximately $k$. Observe, however, that the rule does not
charge a lens much longer than $k$ and charges it very lightly if the lens is much
shorter than $w=w(k)$. This is enough for us to set up the scales of the different phases in such a way
that no crossing is overcharged on average.

We use the following scales for the different phases of our charging: $k=80\alpha^2n/2^{3^u}$, where $u$ is an
integer satisfying $\log\log n/5<u\le\log\log n/2$. We have $M=\lfloor\log\log n/2\rfloor-\lfloor\log\log
n/5\rfloor$ phases.

\begin{lemma}\label{out2} For any crossing point $x\in X'$ the charge leaving
  $x$ by the third rule averaged over the $M$ phases is at most $2$.
\end{lemma}

\proof Each crossing point $x\in X$ is the starting point of exactly two
lenses (no lens starts at a touching point $t\in T$). We bound the average charge sent by a fixed lens $a^*$ by $1$. Let
$l$ be the length of $a^*$ and let $k_0$ be the smallest scale of a phase
where the lens $a^*$ is charged. We have
$l\le3\alpha^3k_0$ and the total charge $a^*$ sends in this phase is less than $v$. For phases with scale
$k>k_0$ we bound the charge $a^*$ sends by
$vl/w(k)\le3v\alpha^3k_0/w(k)$.
Here $w(k)\ge3\alpha^3k_0$ and the value of $w(k)$ exponentially grows as we
consider larger scales. Thus the total
charge $a^*$ sends in all the phases is at most $3v$. With our choice of the
parameters $M\ge3v$ and this proves the estimate claimed. \qed

\vspace{-0.3cm}
\subsection{Total charge received}

We start with a simple observation that will allow us to speak about ``the
next $k$ sad points'' after a poor point on a curve:

\begin{claim}\label{gek}
If a curve $a\in A$ contains at most $k$ sad points, then none of them is poor.
\end{claim}

\proof Clearly, if $|a\cap T'|\le k$, then every intersection point in $X'\cap
a$ sends a charge of $1/k$ to every point in $T'\cap a$ according to the first
rule. The claim follows as there are at least $n-1$ intersection points on
$a$. \qed

Our goal in this subsection is to prove the following lemma.

\begin{lemma}\label{in}
Every sad touching point receives a total charge of at least $\alpha$ in
every phase. 
\end{lemma}

\proof For the proof of this lemma we fix the phase with scale $k$ and we also
fix a single sad  touching point $t\in T'$. We assume for contradiction that
$t$ receives a total charge of less than $\alpha$.
Note first that our assumption implies that $t$ is poor.

\begin{figure}[htbp]
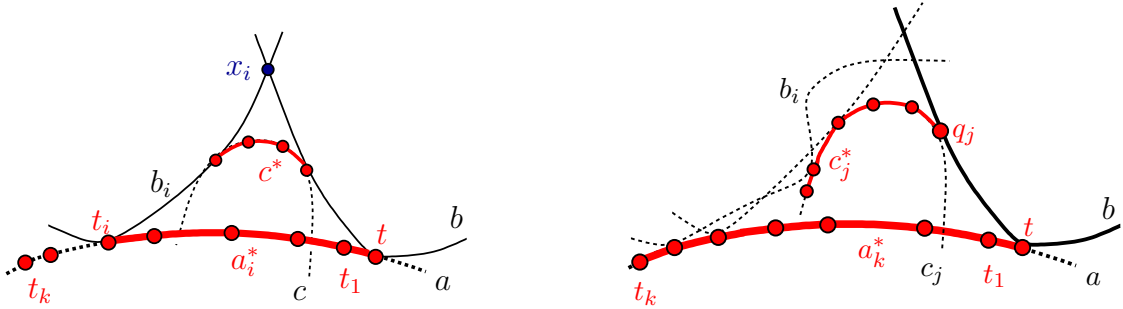

\begin{center}
\input{FastArc.pstex_t}\hspace{2cm}\input{GoodPoint.pstex_t}
\caption{\small Left: A poor arc $c^*$ is $i$-fast if it has apex $x_i$ and starts on the arc of $b$ from $x_i$ to $t$. Right: A good point $q_j$ together with the adjacent arc $c_j^*\subset c_j$. Notice that $c_j$ must meet $a$ within $a^*_k$, and $c_j^*$ contains all the fast arcs that start at $q_j$.}
\label{Fig:GoodPoint}
\end{center}
\end{figure}

Let $t$ be the point
where the curves $a,b\in A$ touch. Let $t_1,t_2,\dots,t_k$ be the first
$k$ sad touching points after $t$ along $a$. By Claim~\ref{gek} they exist. For $1\le i\le k$ let $a_i^*$ be the arc of $a$ from $t$ to $t_i$, let
$x_i$ be the apex of $a_i^*$ and let $b_i$ be the curve in $A$ that touches $a$ at $t_i$. We call a poor arc
{\em$i$-fast} if it starts at a point in the arc from $x_i$ to $t$ and has $x_i$ as its apex, see Figure \ref{Fig:GoodPoint} (left).  We call an arc
{\em fast} if it is $i$-fast for some $1\le i\le k$.

We continue the proof of Lemma~\ref{in} through a series of small claims.

\begin{claim}\label{fast} There are more than $k^2/(2\alpha)$ fast arcs.
\end{claim}

\proof Note that $a_i^*$ itself is $i$-fast and therefore $t$ receives a charge of $2\alpha/k$ from
$x_i$ according to
the third charging rule unless there are more than $k/\alpha$ $i$-fast arcs. As $t$ receives a total charge
of less than $\alpha$ we must have more than $k/\alpha$ $i$-fast arcs for each of more than $k/2$
different values of $i$. This proves the claim. \qed

Note that all fast arcs start at a sad touching point on $b$ and $k$ of them start at $t$. We call a point
{\em good} if at least $k/(4\alpha^2)$ fast arcs start there. Let us name the good points $q_1,\dots,q_L$
in the order they appear on $b$ starting at $q_1=t$ and going along $b$ in reverse direction. For
$1\le j\le L$ let $c_j\in A$ be the curve that touches $b$ at $q_j$ and let $c_j^*$ be the unique arc
from $q_j$ to a point in $c_j\cap T'$ of length exactly $k+1$. The existence
follows from Claim~\ref{gek}. In particular, we have $c_1^*=a_k^*$. See Figure \ref{Fig:GoodPoint} (right).

\begin{claim}\label{goodpoor}
With the previous notation, 
\begin{enumerate}
\item All good points are poor.

\item We have $|c_j^*\cap X'|<\alpha k$ for all $1\le j\le L$.

\item The number of good points is $L\le\alpha k$.

\item At least $k^2/(4\alpha)$ fast arcs start at a good point.

\item Any $i$-fast arc that starts at one of the good points $q_j$ ends at the point where $b_i$ touches
$c_j$ and it is contained in $c_j^*$.
\end{enumerate}
\end{claim}

\noindent{\em Proof:} The first statement holds as any fast arc starts at a poor point by definition.

The second statement follows as each point in $c_j^*\cap X'$ sends a charge of
$1/k$ to the poor point $q_j$ by the first rule.

We prove the third statement in a stronger form: the same bound holds for the number $L'$ of all
starting points of fast arcs. Let $c^*$ be an $i$-fast arc starting at $q\ne t$. The curve $c\in A$ containing
$c^*$ must intersect $a$ and therefore it must escape the triangle like region
bounded by the arc of $b$
from $x_i$ to $t$, $a_i^*$ and the arc of $b_i$ from $t_i$ to $x_i$. It cannot
cross $b$ or $b_i$, so it must leave through (or touch) $a_k^*$ ``using up''
at least one of the at most $\alpha k$ intersection points on
$a_k^*$. Therefore we have $L\le L'\le\alpha k$.

To see the fourth statement note that there are at least $k^2/(2\alpha)$ fast
arcs by Claim~\ref{fast}, but less than $L'k/(4\alpha^2)\le k^2/(4\alpha)$
fast arcs start in points that are not good.

For the final statement note that an $i$-fast arc has length at most $k+1$ by
definition, so if it starts at $q_j$, then it must be contained in
$c_j^*$. The curve $b_i$ must touch $c_j$ at the end point of this arc because
the apex of the arc is on $x_i$.
\qed

We call an arc $z^*$ {\em short} if $|z^*\cap T|\le2\alpha k$. Note that while
the length counts sad touching points on an arc in this definition we count
all touching points. We say that the
good points $q$ and $q'$ are {\em close}, if either the arc
of $b$ from $q$ to $q'$ or the arc from $q'$ to $q$ is short.

\begin{claim}\label{notclose}
Let $q$ and $q'$ be good points. If the arc $b^*$ from $q'$ to $q$ is short,
then $|b^*\cap X'|\le2\alpha(\alpha_1+1)k$.

If $q_j$ and $q_{j'}$ are not close, then the arcs $c_j^*$ and $c_{j'}^*$ are
disjoint.
\end{claim}

\proof The first claim holds simply because $q$ is a sad touching point, so
the arc $b^*$ ending there must have crossing-to-touching ratio below
$\alpha_1$.

For the second claim assume $c_j^*$ and $c_{j'}^*$ intersect and let $W$ be a
Jordan curve connecting $q_j$ to $q_{j'}$ along part of $c_j^*$ and
$c_{j'}^*$. Consider the two arcs $b^*$ and $b'^*$ that $b$ is cut by $q_j$
and $q_{j'}$. By our assumption neither of these arcs is short, so each has
more than $2\alpha k$ distinct curves touching it. As $W\cap X'\le2\alpha k$ we
must have a curve $z\in A$ touching $b^*$ that is disjoint from $W$ and
similarly we have another curve $z'\in A$ touching $b'^*$ and also
disjoint from $W$. Now $b$ and $W$ separate $z$ 
and $z'$ contradicting that they (as any two curves in $A$) must
intersect. The contradiction finishes the proof of the claim. \qed

We call the distinct good points $q_j$ and $q_{j'}$ {\em mingled} if
$|c_j^*\cap c_{j'}^*|>\alpha^2k/w$.

\begin{claim}\label{mingled}
Mingled points are close. A good point is mingled with at most $w/\alpha$
other good points.
\end{claim}

\proof The first statement follows directly from Claim~\ref{notclose}. The
second statement follows from the statement of Claim~\ref{goodpoor} that
$c_j^*$ contains at most $\alpha k$ crossing points in total. \qed

For a good point $q$ let $I_q$ stand for the set of indices $1\le i\le k$ with
an $i$-fast arc starting at $q$. Similarly, for $1\le i\le k$ let $Q_i$ stand
for the set of good points $q$ where an $i$-fast arc starts.

\begin{claim}\label{close} Let $1\le j<j'\le L$ be such that the arc from
$q_{j'}$ to $q_j$ is short, but $q_j$ and $q_{j'}$ are not mingled. Then
$|I_{q_j}\cap I_{q_{j'}}|<6\alpha^2k/\sqrt v$.
\end{claim}

\proof For simplicity we write $q$ and $q'$ for $q_j$ and $q_{j'}$,
respectively. Similarly, we write $c$, $c'$, $c^*$ and $c'^*$ for $c_j$,
$c_{j'}$, $c_j^*$ and $c_{j'}^*$, respectively. We write $b^*$ for the short arc
from $q'$ to $q$.

\begin{figure}[htbp]
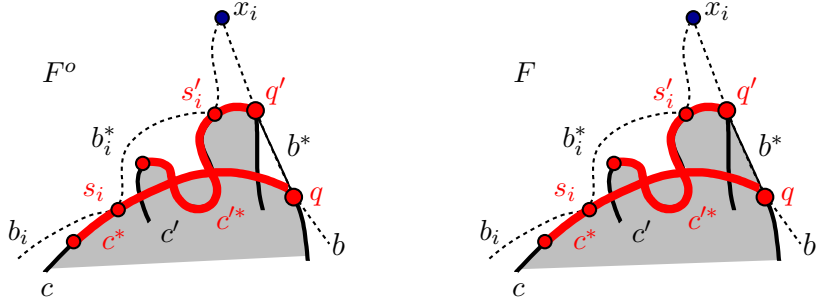

\begin{center}
\input{Faces.pstex_t}\hspace{2cm}\input{Faces1.pstex_t}
\caption{\small Proof of Claim \ref{close}. The curve $b$ lies in the single face $F^o$ of the arrangement of $c$ and $c'$.
The cell $F^o$ contains all the curves $b_i$ with $i\in I$ (left). The arc $b^*$ from
$q'$ to $q$ splits $F^o$ into two sub-faces, with all the touching points
$s_i$ and $s'_i$ lying on the boundary the sub-face $F\subset F^o$ to the left
of $b^*$ (right).}
\label{Fig:Faces}
\end{center}
\end{figure}

Refer to Figure \ref{Fig:Faces}.
Consider the arrangement of the curves $c$ and $c'$. The curve $b$
touches both of these curves so it must be contained in a single face $F^o$ of
the arrangement and this face is to the right of both $c$ and $c'$. The
boundary of $F^o$ is a simple closed Jordan curve which we denote by
$W^o$. Clearly, $W^o$ consists of alternating arcs of $c$ and $c'$ each
consistently oriented with $F^o$ to the right of them. The arcs of $c\cap W^o$ appear in the same cyclic order along $c$ and $W$, and a similar statement
is true for the segments of $c'\cap W^o$. Note, however, that outside $F^o$
the curves $c$ and $c'$ can behave wildly and all sorts of extra intersections
can occur even between $c^*$ and $c'^*$.

The arc $b^*$ splits this face in two, let $F$ stand for the
side of $F^o$ containing $b\setminus b^*$ (that is, on the left from
$b^*$). Let $W^1$ stand for the part of $W^o$ on the boundary of $F$. Clearly,
this is a Jordan curve connecting $q'$ to $q$ and all its segments coming
from $c$ and $c'$ are consistently oriented from $q'$ to $q$.

Let us write $I=I_q\cap I_{q'}$. In what follows, we can assume that $I$ is not empty.
For each $i\in I$ the curve $b_i$ touches both $c$ and $c'$ and intersects
$b$, so it is confined to the face $F^o$. Let $s_i\in c^*$ and
$s'_i\in c'^*$ be the points where $b_i$ touches the boundary of $F^o$. Recall
that $x_i$ is an intersection point of $b_i$ and $b$. Let $b_i^*$ be the arc
of $b_i$ from $s_i$ to $s'_i$ or vice versa,
whichever does not contain $x_i$.
Note that by the definition of the apex $b_i^*$ and $b$ do not intersect and
$x_i$ is the first point on $b_i$ after $b_i^*$ that belongs
to $b$. Furthermore, we must have $x_i\notin b^*$. This shows that the interior
of $b_i^*$ is inside the face $F$ and the endpoints $s_i$ and $s'_i$ are on the
boundary of $F$, and, consequently, also on $W^1$.

As $c'^*$
starts at $q'$ we have a segment $U'$ of $W^1$ starting at $q'$ with $c'^*\cap
W^1=c'\cap U'$. Although $c^*$ starts at $q$ we similarly have $c^*\cap
W^1=c\cap U$ for an interval $U$ of $W^1$ starting at $q'$. Indeed, $c^*\cap W^o$ is
confined to an interval of $W^o$ starting at $q$ and we can choose $U$ to be
the intersection of this interval with $W^1$. 

\begin{figure}[htbp]
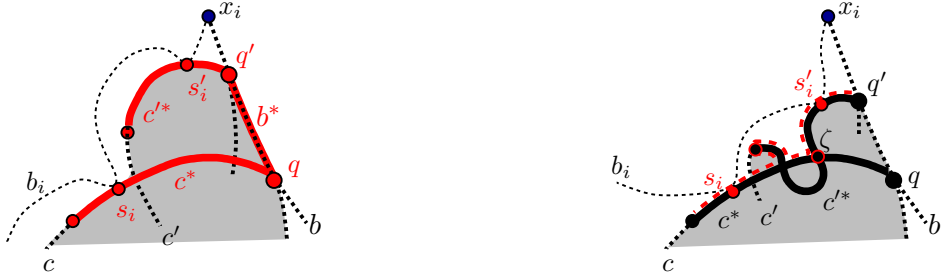

\begin{center}
\input{TraceDisjointFix1.pstex_t}\hspace{4cm}\input{TraceIntersectingFix.pstex_t}
\caption{\small Proof of Claim \ref{close}: Constructing the curve $W$ which
  traces all the portions of $c^*$ and $c'^*$ that appear on the boundary of the face $F$. Left: If $c^*$ and $c'^*$ are disjoint, then $W$ includes the arc $b^*$ from $q'$ to $q$. Right: The curves $c^*$ and $c'^*$ intersect. Notice that the first intersection $\zeta$ of $c'^*$ with $c$ must belong to $c^*$.}  
\label{Fig:TraceW}
\end{center}
\end{figure}

We define a simple open Jordan curve $W$ that contains 
$(c^*\cup c'^*)\cap W^1$ and consists only of segments of $c^*$,
$c'^*$ and possibly one additional segment. In particular, $W$ contains all the points $s_i=c^*\cap b_i$ and $s'_i=c'^*\cap b_i$ for $i\in I$. Refer to Figure \ref{Fig:TraceW}.

In case $c'^*$ does not intersect $c$, it lies entirely on $W^o$. In particular,  $c^*$ and $c'^*$ are disjoint. Hence, we can simply take $W$ to be
the union of $c^*$, $c'^*$ and $b^*$.

In case $c'^*$ intersects $c$, we consider the first such intersection point $\zeta$ along $c'^*$, at which $c'^*$ leaves $W^o$, and notice that $\zeta$ must belong to $c^*$.
Indeed, assume for a contradiction that $\zeta$ belongs to $c\setminus c^*$.
Since (i) the order of the segments of $c\cap W^o$ along $W^o$ is consistent with their order along $c$, (ii) $W^1$ ends at the starting point $q$ of $c^*$, and (iii) $W^1$ begins at $q'$,
the last appearance of $c^*$ along $W^o$ is also contained in $W^o\setminus W^1$. However, in that case $c^*$ can never show
up on $W^1$, contrary to $I\neq \emptyset$ (and thus the choice of $s_i$ on $W^1\cap c^*$ for $i\in I$).

We can assume, then, that the first intersection $\zeta$ of $c'^*$ with $c$ lies on $c^*$, and it is the first appearance of $c$ and $c^*$ along $W^1$.
Our curve  $W$ starts with the shorter of the
two intervals $U$ and $U'$. Suppose this is $U$. Then we have $c^*\cap
W^1\subseteq U\subseteq c^*\cup c'^*$. If the end point $u$ of $U$ is in
$c'^*$ we simply add the part of $c'^*$ after $u$ to obtain $W$. If
$u\notin c'^*$, then we continue $W$ by retracing the last
segment of $U$ very close, but slightly outside $F$ till we meet $c'^*$ and
then add the remaining part of $c'^*$ to finish $W$. We trace $W$ symmetrically if
the shorter segment is $U'$: We follow the remaining part of $c^*$,
possibly using a reverse segment slightly outside $F$ to connect the two
parts. Note that the only case when this construction does not work is when
$U'$ is disjoint from $c^*$, but in this case $c'^*$ is contained in $W^1$ and
is disjoint from $c^*$ and we handled this case separately. Note also, that as
the part of $W$ outside (the very small neighborhood of) $F$ is contained by
just one of $c^*$ or $c'^*$, the curve $W$ is simple: no self-intersections
occur.
\smallskip

Let us define $B$ as the set of points where $W$ is touched or crossed by a
curve in $A$. We have $X'\cap W\subseteq B$ but $B$ may contain further
crossing points along its reverse segment, if such a segment exists. Still, the
overall size of $B$ is at most $\le2\alpha(\alpha_1+1)k+2\alpha
k=2\alpha^2k$ by Claims~\ref{notclose} and \ref{goodpoor}. We
will apply Lemma~\ref{Lemma:RedBlue} to $B$ with a non-uniform weight function
$w_0$. We set $w_0(x)=w+1$ for $x\in B\cap c^*\cap c'^*$ or (in the disjoint
case) if $x=q$. We set $w_0(x)=1$ otherwise. For the total weight we have
$w_0(B)\le3\alpha^2k$ as $q$ and $q'$ are not mingled.

Let $W_i$ be the portion of $W$ between the touching points $s_i$ and $s'_i$
and let us write $l_i=w_0(B\cap W_i)$. Let $R_0=\{s_i\mid i\in I\}$ and set
$\lambda=3\sqrt v$. Let us write
$R=\{s_i\mid i\in I, |W_i\cap R_0|\le l_i/\lambda\}$ and $\hat I=\{i\in I\mid
s_i\notin R\}$. We have $|I|=|R|+|\hat I|$. In what follows we bound $|R|$ and
$|\hat I|$ separately.

To bound $|R|$ we apply Lemma~\ref{Lemma:RedBlue} for the curve $W$, the sets
$R$ and $B$, the weight function $w_0$ and the parameter $\lambda$. The
condition is satisfied as for the interval $W_i$ ending at
$s_i\in R$ we have $|W_i\cap R|\le|W_i\cap R_0|\le l_i/\lambda=w_0(W_i\cap
B)/\lambda$. From the Lemma we conclude that
$$|R|\le\frac{3w_0(B)}\lambda\le\frac{9\alpha^2k}\lambda.\eqno{(2)}$$

\begin{figure}[htbp]
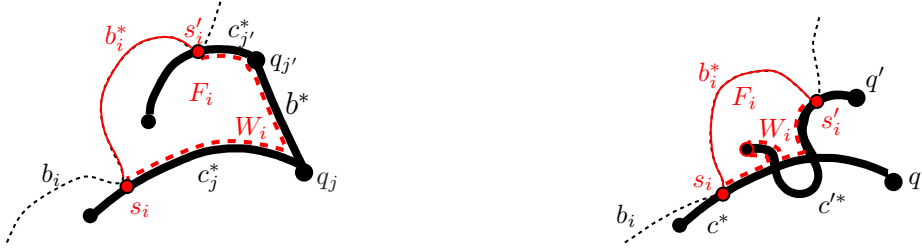

\begin{center}
\input{TraceDisjoint2Fix1.pstex_t}\hspace{4cm}\input{TraceIntersecting2Fix.pstex_t}
\caption{\small The portion $W_i$ of $W$ between $s'_i=b_{i}\cap c'$ and $s_i=b_i\cap c$ is traced. 
The face $F_i$ of $\reals^2\setminus (b_i\cup W_i)$ lies to the right of $b_i$. Left: The arcs $c_j^*$ and $c_{j'}^*$ are disjoint so both $W$ and $W_i$ must include $b^*$.  Right: The scenario where $c_j^*$ and $c_{j'}^*$ intersect.}
\label{Fig:TraceW1}
\end{center}
\end{figure}

For $i\in I$ consider $W_i\cup b_i^*$ and let $F_i$ be the
side of this closed Jordan curve to the right of $b_i^*$; see Figure \ref{Fig:TraceW1}. Note that $b$ is
disjoint from the interior of $F_i$.

We claim that $b_i^*$ contains at most 
$l_i-w$ points of $T$, and at most $\alpha_1(l_i-w)$ points of $X$.
The bound on $|b_i^*\cap T|$ follows from the fact that any curve in $A$
touching $b_i^*$ must intersect $W_i$. Indeed, any such curve is in $F_i$ in a
small neighborhood around the point where it touches $b_i^*$, but as $b$ is
disjoint from the interior of $F_i$ it must intersect the boundary of
$F_i$ at a point other than the touching with $b_i^*$. This intersection must
be in $B\cap W_i$. We have $|B\cap W_i|\le l_i-w$ as $W_i$ contains at
least one of the heavy points with weight $w+1$. The bound on $|b_i^*\cap X|$
follows since the endpoint $s_i$ of $b_i^*$ is sad.
\medskip

\begin{figure}[htbp]
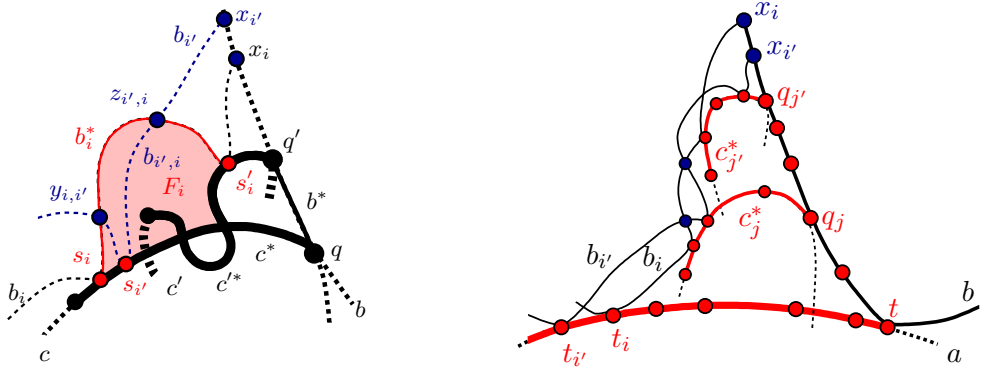

\begin{center}
\input{TraceIntersecting3Fix.pstex_t}\hspace{2cm}\input{TowerFarPoints.pstex_t}
\caption{\small Left: The segment $W_i$ contains a touching $s_{i'}$ between $c$ with $b_{i'}$, for $i'\in I\setminus \{i\}$. Each such $b_{i'}$ yields a lens $b_{i',i}$ of length at most $\alpha(l_i-w)$. Right: Proof of Claim \ref{tower}.
The good points $q_j,q_{j'}\in S$ are not close, and their arcs $c^*_j$ and $c^*_{j'}$ form a pair of neighboring ``teeth''  in the comb $\Gamma$.}
\label{Fig:TraceLens}
\end{center}
\end{figure}

Let us consider $i,i'\in I$ with $i\ne i'$ and $s_{i'}\in W_i$.
Follow $b_{i'}$ from $s_{i'}$ in both directions. It starts out inside $F_i$
and eventually has to reach $b$ that is disjoint of the interior of $F_i$. As
the part of $b_{i'}$ around $s_{i'}$ is in $F$, the first intersection with
$b$ in either direction is outside $b^*$. Since it is not enough for $b_{i'}$
to just touch $b^*$, 
$c_{i'}$ has to cross the boundary of $F_i$. This first crossing point in
either direction must be on $b_i^*$, since it cannot be on $b^*$ and $b_{i'}$
touches both $c_j$ and $c_{j'}$. Let us call these crossing points $y_{i,i'}$
and $z_{i',i}$ with the arc $b_{i',i}^*$ from $y_{i',i}$ to $z_{i',i}$ along
$b_{i'}$ being inside $F_i$ and containing $s_{i'}$. (See Figure \ref{Fig:TraceLens} (left).)

Note that $b_{i',i}^*$ is a lens. The length $l_{i',i}$ of this lens is at most
$\alpha(l_i-w)$. Indeed, any curve touching $b_{i,i'}^*$
must intersect the boundary of $F_i$ in order to intersect $b$, but the
total number of points on this boundary, where a curve of $A$ intersects it is
at most $|B\cap W_i|+|b_i^*\cap X'|\le(\alpha_1+2)(l_i-w)=\alpha(l_i-w)$.

Each such lens $b_{i',i}$ sends a
charge of $v/((l_{i',i}+w)k)$ to $q$ by the second
charging rule. Indeed, this rule applies to $b_{i',i}$ as its length is
$l_{i',i}\le\alpha l_i\le\alpha|B|\le 3\alpha^3k$ and
the arc from $q$ to $s_{i'}$ satisfies the requirements. The amount of the
charge sent is
$$\frac v{(l_{i',i}+w)k}\ge\frac v{(\alpha(l_i-w)+w)k}\ge\frac
v{\alpha(l_i-\lambda)k}.$$

If we further assume that $i\in\hat I$, then we have more than $l_i/\lambda$
choices of $i'\in I$ with $s_{i'}\in W_i$. One of these choices is $i'=i$, but
more than $l_i/\lambda-1$ other choices will give rise to lenses $b_{i',i}$,
each sending a charge of at least $v/(\alpha(l_i-\lambda)k)$ to
$q$. The total of these charges for a fixed $i\in\hat I$ is at least
$v/(\alpha\lambda k)$ and for all $i\in\hat I$ it is at least
$|\hat I|v/(\alpha\lambda k)$.

We know that $q_j$ is poor, so this charge does not reach the threshold of
$\alpha$. As a consequence, we have
$\displaystyle |\hat I|\le\alpha^2\lambda k/v.$
To finish the proof of this claim we use this last estimate, Equation~(2),
the fact $|I|=|R|+|\hat I|$ and substitute $\lambda=3\sqrt v$. \qed
\smallskip

\begin{claim}\label{fewclose}
Let $1\le j<j'\le L$ be such that the arc $b^*$ from $q_{j'}$ to $q_j$ is
short. The number of good points in $b^*$ is at most $50\alpha w$.
\end{claim}

\proof Let $S$ be the set of good points in $b^*$. We have
$\displaystyle |S|\frac k{4\alpha^2}\le\sum_{q\in S}|I_q|=\sum_{i=1}^k|Q_i\cap S|$
and by the Cauchy--Schwarz inequality
$\displaystyle \frac{k|S|^2}{16\alpha^4}\le\sum_{i=1}^k|Q_i\cap S|^2=\sum_{q,q'\in
S}|I_q\cap I_{q'}|.$
We use the trivial bound $|I_q\cap I_{q'}|\le k$ if $q=q'$ or if $q$ and $q'$ are
mingled. By Claim~\ref{mingled} there are at most $(w/\alpha+1)|S|$ such
terms. The remaining terms can be bounded by $6\alpha^2k/\sqrt v$ using
Claim~\ref{close}. We obtain
$\displaystyle \frac{k|S|^2}{16\alpha^2}\le k(w/\alpha+1)|S|+\frac{6|S|^2\alpha^2k}{\sqrt v}.$
Substitute $v=21000\alpha^8$ and the claim follows. \qed.

\begin{claim} \label{tower}
For distinct indices $1\le i,i'\le k$ the curves $b_i$ and $b_{i'}$ have at least
$|Q_i\cap Q_{i'}|/(100\alpha w)-1$ crossing points.
\end{claim}

\proof Let $Q=Q_i\cap Q_{i'}$. We select a subset $S\subseteq Q$ with no two
close points greedily: we consider the elements $q\in Q$ along $b$ in reverse
direction starting
from $t$ and include them in $S$ unless $q$ is close to a point already in
$S$. No short arc of $b$ avoiding $t$ contains  more than $50\alpha w$
good points by Claim~\ref{fewclose}, thus we have
$|S|\ge|Q|/(50\alpha w)-1$. Let $b^*$ be the
arc of $b$ from the point we put in $S$ last to the point we put there
first. By Claim~\ref{notclose} the arcs $c_j^*$
corresponding to the points $q_j\in S$ are pairwise disjoint. Let $\Gamma$ be
the comb like arrangement of these arcs together with $b^*$, see Figure \ref{Fig:TraceLens} (right). Let
$b_i^*$ be the maximal arc on $b_i$ from a touching point where an $i$-fast
arc ends to the apex $x_i$ of the $i$-fast arcs. Clearly, $b_i^*$ touches all
the ``teeth'' of the comb $\Gamma$, but it does not intersect its spine
$b^*$. This implies that $b_i^*$ touches the teeth in the same order as the
$b^*$. This is also true for the analogously defined arc $b_{i'}^*$ of
$b_{i'}$. Consider two neighboring teeth of the comb $\Gamma$. Clearly, either the
part of $b_i^*$ between the corresponding touching points is crossed by
$b_{i'}$ or the part of $b_{i'}^*$ between the corresponding touching points
is crossed by $b_i$. As the segments of $b_i^*$ between touchings of
consecutive teeth are disjoint this represent at least $(|S|-1)/2$ crossings
between $b_i$ and $b_{i'}$.
\qed.

\medskip
Having proved these claims we return to the proof of Lemma~\ref{in}. We started the proof by
assuming the lemma fails, so we need to arrive to a contradiction to finish the proof.

By Claim~\ref{goodpoor} we have many fast arcs starting at good points, namely
$ \sum_{j=1}^L|I_{q_j}|\ge\frac{k^2}{4\alpha}.$
Using $L\le\alpha k$ (Claim~\ref{goodpoor} again) and the Cauchy--Schwarz inequality we get
$ \frac{k^3}{16\alpha^3}\le\sum_{j=1}^L|I_{q_j}|^2=\sum_{1\le i,i'\le
k}|Q_i\cap Q_{i'}|.$
Subtracting the contribution of the $i=i'$ case and dividing by $2$ we get
$\sum_{1\le i<i'\le k}|Q_i\cap Q_{i'}|\ge\frac{k^3}{40\alpha^3}.$
Claim~\ref{tower} shows that this lower bound on the left hand side provides a lower bound on the
number of the crossing points between the curves $b_i$. We find that at least
$k^3/(4000\alpha^4w)-k^2/2$ such crossing points exist. With the prior choice of parameters, 
this contradicts the assumption that the total number of crossings
satisfies $|X|<\alpha_1{n\choose2}$. The contradiction proves Lemma 5. \qed

\vspace{-0.3cm}
\subsection{Wrap up.} Finishing the proof of Theorem~\ref{main} is simple once we have
Lemmas~\ref{happy}, \ref{out13}, \ref{out2}
and \ref{in}. Considering
all the charges in all the $M$ phases of our scheme every sad touching point $t\in T'$ receives a charge of at
least $c_{\mbox{\scriptsize in}}=\alpha M$ by Lemma~\ref{in}. For an intersection
point $x\in X'$ the total charge sent
out is at most $c_{\mbox{\scriptsize out}}=10M$ by Lemmas~\ref{out13} and \ref{out2}. Comparing the total charges sent and received we obtain
$$\frac{|X'|}{|T'|}\ge\frac{c_{\mbox{\scriptsize in}}}{c_{\mbox{\scriptsize
out}}}=\alpha/10.$$
We have $|T'|\le 10|X'|/\alpha$ from the line above and
$|T|-|T'|\le6|X|/\alpha_1$ from Lemma~\ref{happy}. In total we have
$|T|\le16|X'|/\alpha_1$ and the statement of the Theorem~\ref{main} follows.

\end{document}